\documentclass[11pt]{article}

\usepackage[latin1]{inputenc}
\usepackage{graphicx}
\usepackage{amsmath}
\usepackage{amssymb}
\usepackage{tikz}
\usepackage{caption}
\usepackage{hyperref}
\usepackage{eurosym}
\usepackage{algorithm}
\usepackage{cite}
\usepackage{hhline}
\usepackage{commandes}
\usepackage{jlcode}

\newcommand{\stamp}{t}

\newcommand{\pv}{g}
\newcommand{\demand}{d}
\newcommand{\soc}{x}
\newcommand{\ConstraintSet}{\mathcal{U}}
\newcommand{\StageCost}{L}
\newcommand{\TotalCost}{\mathcal{L}}

\newcommand{\hist}{h}
\newcommand{\HIST}{\mathbb{H}}

\newcommand{\ExchangeEnergy}{e}

\newcommand{\SimulationData}{\mathscr{S}}
\newcommand{\CalibrationData}{\mathscr{C}}
\newcommand{\BaseData}{\mathscr{D}}

\newcommand{\batt}{(c, \overline{l}, \rho_c, \rho_d)}

\newcommand{\emsx}{\texttt{EMSx.jl}}

\newcommand{\site}{i}
\newcommand{\SITE}{I}

\newcommand{\uspan}{\na{0, \ldots, T{-}1}}
\newcommand{\wspan}{\na{1, \ldots, T}}

\newcommand{\Gain}{G}
\newcommand{\Score}{\mathcal{G}}
\newcommand{\cardinal}[1]{\module{#1}}

\newcommand{\MPC}{\mathrm{\scriptscriptstyle MPC}}
\newcommand{\OLFC}{\mathrm{\scriptscriptstyle OLFC}}

\newcommand{\SDPt}{\textsc{\tiny SDP}}
\newcommand{\SDP}{\textsc{SDP}\,}
\newcommand{\SDPARt}{\textsc{\tiny SDP-AR}}
\newcommand{\SDPAR}{\textsc{SDP-AR}\,}

\renewcommand{\scenario}{\sigma}
\renewcommand{\SCENARIO}{\mathbb{S}}
\newcommand{\offline}{\mathrm{\scriptscriptstyle off}}
\newcommand{\online}{\mathrm{\scriptscriptstyle on}}

\graphicspath{{Figures/}}
\usepackage{a4wide}
\def\keywords#1{\par\addvspace\medskipamount{\rightskip=0pt plus1cm
\def\and{\ifhmode\unskip\nobreak\fi\ $\cdot$
}\noindent{\bfseries Keywords}\enspace\ignorespaces#1\par}}
\def\myacknowledgement{\par\addvspace{17pt}\small\rmfamily
\trivlist\if!{Acknowledgements}!\item[]\else
\item[\hskip\labelsep
{\bfseries{Acknowledgements}}]\fi}

\newenvironment{acknowledgements}{\begin{myacknowledgement}}
{\end{myacknowledgement}}

\title{EMSx: a numerical benchmark \\ for energy management systems}

\author{Adrien~Le~Franc\thanks{CERMICS, \'{E}cole des Ponts ParisTech, France.}
  \thanks{Efficacity}
  \and
  Pierre~Carpentier\thanks{UMA ENSTA Paris} \and
  Jean-Philippe~Chancelier\footnotemark[1]\and
  Michel~De~Lara\footnotemark[1]
}

\begin{document}

\maketitle

\begin{abstract}
	Inserting renewable energy in the electric grid in a decentralized manner
	is a key challenge of the energy transition.
	However, at local scale, both production and demand display erratic behavior,
	which makes it \noir{challenging} to match them. 
	It is the goal of Energy Management Systems (EMS) to achieve such balance at
	least cost.
	We present EMSx, a numerical benchmark for testing control algorithms for the
	management of electric microgrids equipped with a photovoltaic unit and an energy
	storage system. 
	EMSx is made of three key components: 
	the EMSx dataset, provided by Schneider Electric, contains a diverse pool of realistic
	microgrids with a  rich collection of historical observations and forecasts;
	the EMSx mathematical framework is an explicit description of the assessment
	of electric microgrid control techniques and algorithms;
	the EMSx software  \emsx\ is a package, implemented in the Julia language, 
	which enables to easily implement a microgrid controller and to test it.
	All components of the benchmark are publicly
	available, so that other researchers willing to test controllers on EMSx may
	reproduce experiments easily.
	Eventually, we showcase the results of standard microgrid
	control methods, including Model Predictive Control, Open Loop Feedback Control
	and Stochastic Dynamic Programming.
      \end{abstract}



  \section{Introduction}
\label{intro}

Inserting renewable energy in the electric grid is a key challenge of the energy
transition. As renewable power units are often coupled with a storage system and envisaged at 
a local scale (where the demand is more erratic than at a global scale),
this makes the management of such electric microgrids delicate.
In this paper, 
we concentrate on Energy Management Systems (EMS), a high level layer of hierarchical control, 
responsible for operating electric microgrids while optimizing security and economic criteria
\cite{olivares2014trends}. 

The design of an EMS for the optimal 
management of a storage system supporting renewable energy integration is a
well-known challenge. Most approaches are based on the mathematical theory of multistage 
deterministic and stochastic optimization, and the relevance of a microgrid control
method is usually assessed on a representative validation setup.
A comparative study of EMS design techniques was conducted in
\cite{rigaut2018stochastic}, which includes Model Predictive Control (MPC,
\cite{garcia1989model,bertsekas2005dynamic}), Open Loop Feedback Control (OLFC,
\cite[Vol.~1, \S6.2]{bertsekas1995dynamic}, sometimes referred to as stochastic
MPC) and Stochastic Dynamic Programming (\SDP
\cite{bertsekas12,puterman94}). Although providing a large benchmark of control
methods, \cite{rigaut2018stochastic} assesses management strategies on a very
specific example, the control of a braking energy storage system in a subway
station. 
Another example of highly specific EMS assessment result is found in
\cite{kriett2012optimal}, where the authors apply MPC for operating a very
detailed residential microgrid configuration, mixing thermal and electric energy
storage, and providing models for every electric device accounting for a
shiftable load in a standard Swiss household. 
Due to the specificity of the case studies, it is not clear that the conclusions
of such works can be generalized, and that they can 
help researchers and practitioners willing to deploy an EMS for a new
context, with new data.

Novel microgrid controller techniques are published regularly, making it
even harder to establish an up-to-date benchmark. 
In \cite{hafiz2019real}, the
authors propose an innovative EMS architecture which combines Recurrent Neural
Networks \cite{gers2002learning} with Stochastic Dual Dynamic Programming
\cite{shapiro2011analysis} for managing solar panels coupled with a
battery. They compare their EMS against a heuristic approach on residential
data collected by a research organization which provides free access to its
data for academic purposes. However, their simulation relies on a
commercialized real-time simulator --- which may improve the realism of the
assessment method, but complicates the reproduction of the testbed. 
A similar
microgrid application is considered in \cite{heymann2016stochastic}, where the
authors also operate a battery for photovoltaic power integration in a novel
framework, modeling the electric load by a stochastic differential
equation. They rely on real data --- whose access is not documented --- from an
experimental site in Chile, and implement their simulation on an open source
numerical solver. 
The same initiative was taken in
\cite{haessig2014computing}, where the authors make their implementation of
SDP publicly available, so that their comparison with a heuristic method for
the design of an EMS could be extended to other data and other techniques.
\medskip

It is in that context that we introduce the EMSx benchmark, composed of three
constituents --- a dataset, a mathematical framework, and a simulation software
--- and designed for the purpose of assessing 
electric microgrid controllers on an open, transparent and unified challenge.  

The first component of the EMSx benchmark is 
a new \emph{dataset} reporting a wide range of contrasted 
microgrid situations. For these data, we rely on \noir{samples} collected by the Schneider
Electric (SE) company on 70 industrial sites. SE develops and commercializes
microgrid controllers, and is interested in challenging its current practices
with state-of-the art optimization methods developed in the academic world. In
this context, SE puts together a dataset \noir{which allows benchmarking} 
different energy management methods in various situations. 
This dataset contains photovoltaic and electric load profiles, 
detailing both historical observations and forecasts, which let us
closely reproduce the online information available to an EMS.

The second component of the EMSx benchmark is 
a \emph{mathematical framework} for the assessment
of electric microgrid control techniques and algorithms. 
This framework consists \noir{of} the mathematical description of 
a microgrid architecture --- made of one resource, one battery, one load ---
and of an economic criterion that depends on how the battery is managed.
Indeed, as we focus on photovoltaic units 
integrated in local power networks --- with uncertainties arising both from 
the electric load and from the photovoltaic power generation ---
introducing an energy storage system can help \noir{to reduce} the energy bill. 
The battery is managed by means of a microgrid controller, for which we propose
a precise mathematical definition and  a score to measure its
aggregated performance, on the one hand across the uncertainties 
unique to a site,  on the other hand across the different sites of the dataset. 

The third component of the EMSx benchmark is the EMSx \emph{benchmark software},
that we have developed and implemented as a Julia \cite{bezanson2012julia} package.
This software makes possible the simulation and 
assessment of a large range of controllers on the testbed defined by our dataset
and our mathematical framework. 
\medskip

We believe that our benchmark is well-suited for
assessing the performance of a large class of control techniques, 
and we illustrate our claim by measuring the performance of a selection 
of controllers --- derived from MPC, OLFC, SDP, and from an extended 
state formulation of a plain \SDP controller that models uncertainties 
with an auto-regressive process (\SDPAR, \cite[\S3.1.1]{Shapiro:2014:LSP:2678054} 
and \cite{lohndorf2019modeling}). 
All in all, EMSx stands out by providing open access to our simulation data and
software, permitting to benchmark microgrid controllers on simulation data 
representative of a large panel of industrial microgrid cases, within a clear and
detailed mathematical framework. 
\medskip

This paper is structured as follows. 
We introduce the EMSx benchmark dataset in Sect.~\ref{sec:data}, 
the EMSx benchmark mathematical framework in
Sect.~\ref{sec:mathematical_framework} --- by providing a microgrid simulation model, 
a definition of a microgrid controller and a score --- 
and the EMSx benchmark \emsx\ package in Sect.~\ref{sec:emsx_package}.
In Sect.~\ref{sec:experiments}, we detail two main classes of 
mathematical techniques to design microgrid controllers,
and we provide the numerical results obtained with the EMSx benchmark. 
The Appendix~\ref{Appendix} gathers additional material
on the numerical experiments of Sect.~\ref{sec:experiments}.

\section{The EMSx benchmark dataset}
\label{sec:data}

We present the Schneider Electric (SE) dataset which offers a large collection of
field data for studying the management of electric microgrids.
The dataset is publicly available at URL \url{https://github.com/adrien-le-franc/EMSx.jl}.

In~\S\ref{Generalities_about_the_Schneider_Electric_70_sites}, we present 
generalities about Schneider Electric's 70 sites.
Then, we detail the content of the dataset by focusing on 
historical observations in~\S\ref{Historical_observations}
and on historical forecasts in~\S\ref{Historical_forecasts}.
Finally, in~\S\ref{Illustrating_how_sites_differ_in_terms_of_predictability},
we illustrate how the 70~sites differ in terms of ``predictability''. 

\subsection{Generalities about Schneider Electric's 70 sites}
\label{Generalities_about_the_Schneider_Electric_70_sites}

Indeed, SE has collected a large \noir{database} of load
profiles on real operated microgrids deployed on a various collection of 70
industrial sites, mainly located in Europe and in the United States
(for data privacy reasons, SE does not provide specific information 
on the origin of the sites). 

Each site is documented with parameters and time series data.
We denote the set of sites by~$\SITE$ (hence $\cardinal{\SITE}=70$).
On each site \( \site \in \SITE \), we provide 
battery parameters $(c^{\site}, \overline{l}^{\site}, \rho_c^{\site},
\rho_d^{\site})$
(see the storage dynamics part in~\S\ref{Microgrid_control_simulation} below).
Regarding time series, 
the database contains historical observations and historical forecasts 
of the energy demand and of photovoltaic generation for each site.
For this latter, however, SE has selected a single photovoltaic profile from a site
located in South Central United States, and has then rescaled this profile for 
each site of the database, taking care \noir{to restore} the balance between the energy 
generation and the load profile. This special treatment of photovoltaic
generation is due to the lack of detailed historical meteorological data for
most of the sites, which is an impediment to compute accurate photovoltaic
forecasts.

We sample the continuous time every 15 minutes, giving,
for each site \( \site \in \SITE \), a \emph{time index
	\( t \in \{1, 2, \ldots,\theta^{\site}{-}1, \theta^{\site}\} \)},
up to the \emph{horizon~$\theta^{\site}$}
(with at least one year of historical observations and forecasts per site).
Every time interval $[t, t{+}1[$ corresponds to 15~minutes.

\subsection{Historical observations}
\label{Historical_observations}

\begin{figure}[htpb]
	\centering
	\begin{minipage}[t]{0.45\textwidth}
          \includegraphics[width=\textwidth, height=4cm]{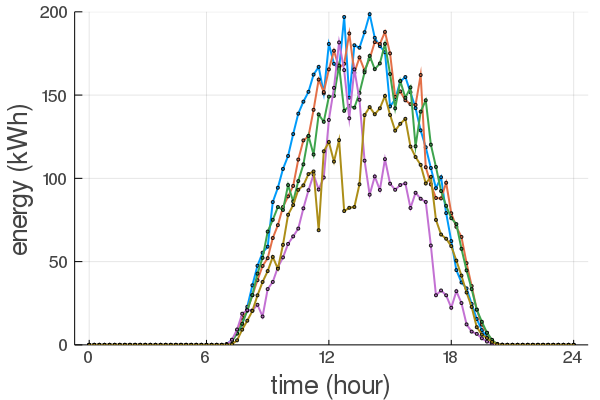}
		\caption{Examples of daily chronicles (historical observations) of photovoltaic generation~\eqref{eq:historical_generation} (from Site~25).}
		\label{fig:pv}
	\end{minipage}
	\hspace{0.3cm}
	\begin{minipage}[t]{0.45\textwidth}
		\includegraphics[width=\textwidth, height=4cm]{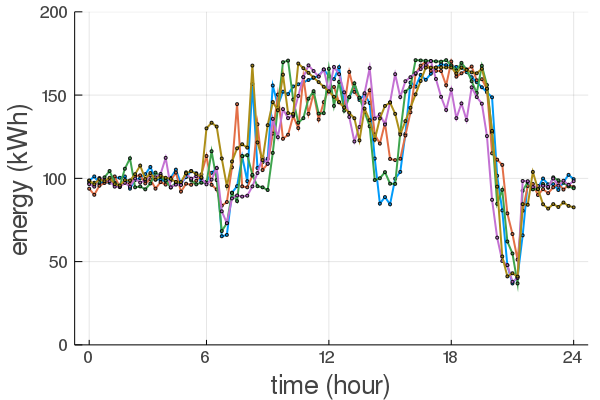}
		\caption{Examples of daily chronicles (historical observations) of energy
			demand~\eqref{eq:historical_demand} (from Site~25)} 
		\label{fig:load}
	\end{minipage}
\end{figure}

We \noir{have} measures of the \emph{photovoltaic generation~$\pv_t$}
and of the \emph{energy demand~$d_t$} over the last 15 minutes,
providing vectors of \emph{historical observations}
for each site \( \site \in \SITE \):
\begin{subequations}
	\begin{align}
	\pv^{\site}
	&=
	\np{ \pv_1^{\site}, \ldots,\pv_{\theta^{\site}}^{\site}}
	\in \RR^{\theta^{\site}} \eqsepv \forall \site \in \SITE
	\eqfinv \label{eq:historical_generation}
	\\
	\demand^{\site}
	&=
	\np{ \demand_1^{\site}, \ldots, \demand_{\theta^{\site}}^{\site}}
	\in \RR^{\theta^{\site}} \eqsepv \forall \site \in \SITE
	\eqfinp \label{eq:historical_demand}
	\end{align}
	\label{eq:historical_observations}
\end{subequations}
We provide examples of observed
daily chronicles for photovoltaic generation in Figure~\ref{fig:pv}, and for
energy demand in Figure~\ref{fig:load}.

\subsection{Historical forecasts}
\label{Historical_forecasts}

On top of observed data, Schneider Electric also provides, every 15 minutes,
\emph{historical forecasts \(\hat{\pv}_{t,t{+}1}^{\site}, \ldots,
	\hat{\pv}_{t,t+96}^{\site}\) of photovoltaic profiles}
and \emph{historical forecasts \(\hat{\demand}_{t,t{+}1}^{\site}, \ldots,
	\hat{\demand}_{t,t+96}^{\site}\) of demand profiles}
for the next 24 hours (hence $96=24 \times 60/15$), hence giving vectors
for all sites \( \site \in \SITE \)
and for all times \( t \in \na{1, \ldots, \theta^{\site}} \)
\begin{subequations}
	\begin{align}
	\hat{\pv}_t^{\site}
	&=
	\np{ \hat{\pv}_{t,t{+}1}^{\site}, \ldots, \hat{\pv}_{t,t+96}^{\site} }
	\in \RR^{96}
	\eqsepv \forall t \in \na{1, \ldots, \theta^{\site}} \eqsepv \forall \site \in \SITE
	\eqfinv
	\\
	\hat{\demand}_t^{\site}
	&=
	\np{ \hat{\demand}_{t,t{+}1}^{\site}, \ldots, \hat{\demand}_{t,t+96}^{\site} }
	\in \RR^{96}
	\eqsepv \forall t \in \na{1, \ldots, \theta^{\site}} \eqsepv \forall \site \in \SITE
	\eqfinv
	\end{align}
\end{subequations}
and sequences of \emph{historical forecasts}
for all sites \( \site \in \SITE \)
\begin{subequations}
	\begin{align}
	\hat{\pv}^{\site}
	&=
	\np{ \hat{\pv}_1^{\site}, \ldots, \hat{\pv}_{\theta^{\site}}^{\site} }
	\in \RR^{96 \times \theta^{\site}}
	\eqsepv \forall \site \in \SITE
	\eqfinv
	\\
	\hat{\demand}^{\site}
	&=
	\np{ \hat{\demand}_1^{\site}, \ldots, \hat{\demand}_{\theta^{\site}}^{\site} }
	\in \RR^{96 \times \theta^{\site}}
	\eqsepv \forall \site \in \SITE
	\eqfinp
	\end{align}
	\label{eq:historical_forecasts}
\end{subequations}
The forecasting method employed by Schneider Electric combines auto-re\-gre\-ssive
models with random forests, and is inspired by the top-level methods used in
GEFCom2014 \cite{hong2016probabilistic}.

By combining the chronicles~\eqref{eq:historical_observations} of \emph{historical observations}
with the chronicles~\eqref{eq:historical_forecasts} of \emph{historical forecasts},
we can thus closely reproduce the information available to an online microgrid controller
operating a real site.

\subsection{Illustrating how sites differ in terms of predictability}
\label{Illustrating_how_sites_differ_in_terms_of_predictability}

The dataset covers a large spectrum of situations regarding variability
and predictability. To assess the predictability of the data at a given site
$\site \in \SITE$, we first introduce 
the \emph{historical net demand observations} 
\begin{subequations}
	\begin{equation}
	z_t^{\site} = \demand_t^{\site} - \pv_t^{\site} \in \RR 
	\eqsepv \forall t \in \{1, \ldots, \theta^{\site}\} 
	\eqsepv \forall \site \in \SITE 
	\eqfinv
	\end{equation}
	deduced from \eqref{eq:historical_observations}, 
	and the \emph{historical net demand forecasts} 
	\begin{equation}
	\hat{z}_t^{\site} = \hat{\demand}_t^{\site} - \hat{\pv}_t^{\site} \in \RR^{96} 
	\eqsepv \forall t \in \{1, \ldots, \theta^{\site}\} 
	\eqsepv \forall \site \in \SITE 
	\eqfinv
	\label{eq:historical_net_demand_forecasts}
	\end{equation}
	deduced from \eqref{eq:historical_forecasts}. 
\end{subequations}
Second, we normalize the historical net demand observations by
\begin{subequations}
	\begin{align}
	&\tilde{z}_t^{\site} = 
	\frac{z_t^{\site} - \underline{z}^{\site}}{\overline{z}^{\site} - \underline{z}^{\site}} 
	\in \nc{0, 1} \eqsepv \forall t \in \na{1, \ldots, \theta^{\site}} 
	\eqsepv \forall \site \in \SITE 
	\eqfinv
	\\
	&
	\mtext{where } \left\{
	\begin{array}{ll}
	\overline{z}^{\site} 
	&= 
	\max_{1\leq t \leq \theta^{\site}}\na{z^{\site}_t} 
	\eqsepv \forall \site \in \SITE 
	\eqsepv 
	\\
	\underline{z}^{\site} 
	&= 
	\min_{1\leq t \leq \theta^{\site}}\na{z^{\site}_t} 
	\eqsepv \forall \site \in \SITE 
	\eqfinv
	\end{array}
	\right.
	\end{align}
\end{subequations}
and we apply the same transformation to the components $\hat{z}_{t,
	t{+}1}^{\site}, \ldots, \hat{z}_{t, t+96}^{\site}$ of the historical net
demand forecasts, yielding
\begin{equation}
\tilde{\hat{z}}^{\site}_{t, t+k} = 
\frac{\hat{z}_{t, t+k}^{\site} - \underline{z}^{\site}}{\overline{z}^{\site} - \underline{z}^{\site}}
\in \RR \eqsepv \forall k \in \na{1, \ldots, 96} 
\eqsepv \forall t \in \na{1, \ldots, \theta^{\site}} 
\eqsepv \forall \site \in \SITE 
\eqfinp
\end{equation}
Thus normalized, predictability can be compared across the pool of sites. 
Third, we define the Root Mean Square Error (RMSE) of the site 
\( \site \in \SITE \) by 
\begin{equation}
\mtext{RMSE}^{\site} = \sqrt{\frac{1}{96 \times \theta^{\site}} 
	\sum_{t=1}^{\theta^{\site}} \sum_{k=1}^{96} 
	\bp{ \tilde{\hat{z}}^{\site}_{t, t+k} - \tilde{z}_{t+k}^{\site} }^2 } 
\eqfinp
\label{eq:RMSE}
\end{equation}

\begin{figure}[htpb]
	\centering
	\includegraphics[width=0.75\textwidth]{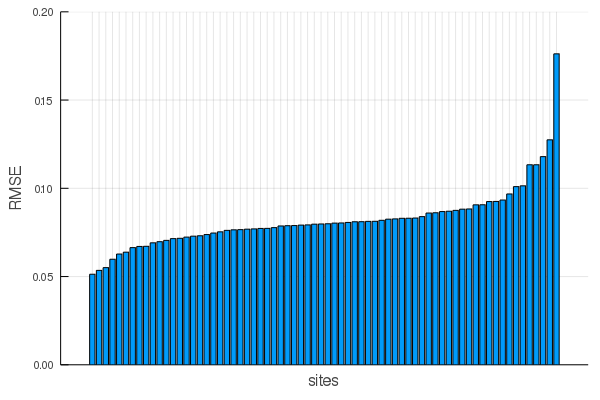}
	\caption{RMSE of the net demand forecast ($Y$-axis) over the pool of 70 sites, ranked in increasing order of forecast error
		along the $X$-axis}
	\label{fig:RMSE}
\end{figure}

The diversity of the forecast error over the pool of 70 sites is \noir{shown}
in Figure~\ref{fig:RMSE}. Here, we ranked sites in increasing order of RMSE. We
observe that the error is quite stable in the wide flat central part of sites
distribution, \noir{except} for the first 20\%
(with low forecast error) and for the last 20\% (with high forecast error) of
the sites. Thus, our dataset represents a diverse pool of microgrids,
offering thus the possibility to assess the performance of microgrid management
techniques on a large range of realistic industrial applications.

\section{The EMSx benchmark mathematical formulation}
\label{sec:mathematical_framework}

After Schneider Electric's dataset presented in Sect.~\ref{sec:data},
the second component of the EMSx benchmark is 
a \emph{mathematical framework} for the assessment
of electric microgrid control techniques and algorithms. 
We formulate a benchmark problem to 
evaluate generic microgrid controllers in~\S\ref{Microgrid_control_simulation}.
Then, we describe the structure of a controller,
and how it is assessed by simulation along a partial chronicle
in~\S\ref{Microgrid_controller}.
In~\S\ref{Designing_and_assessing_a_controller_on_a_given_site},
we detail a score to compare the performance of controllers,
and, in~\S\ref{Assessing_a_controller_design_technique},
we extend it into a score to assess a whole controller-design technique.

\subsection{Microgrid control model}
\label{Microgrid_control_simulation}

\begin{figure}[htpb]
	\centering
	\includegraphics[width=0.6\textwidth]{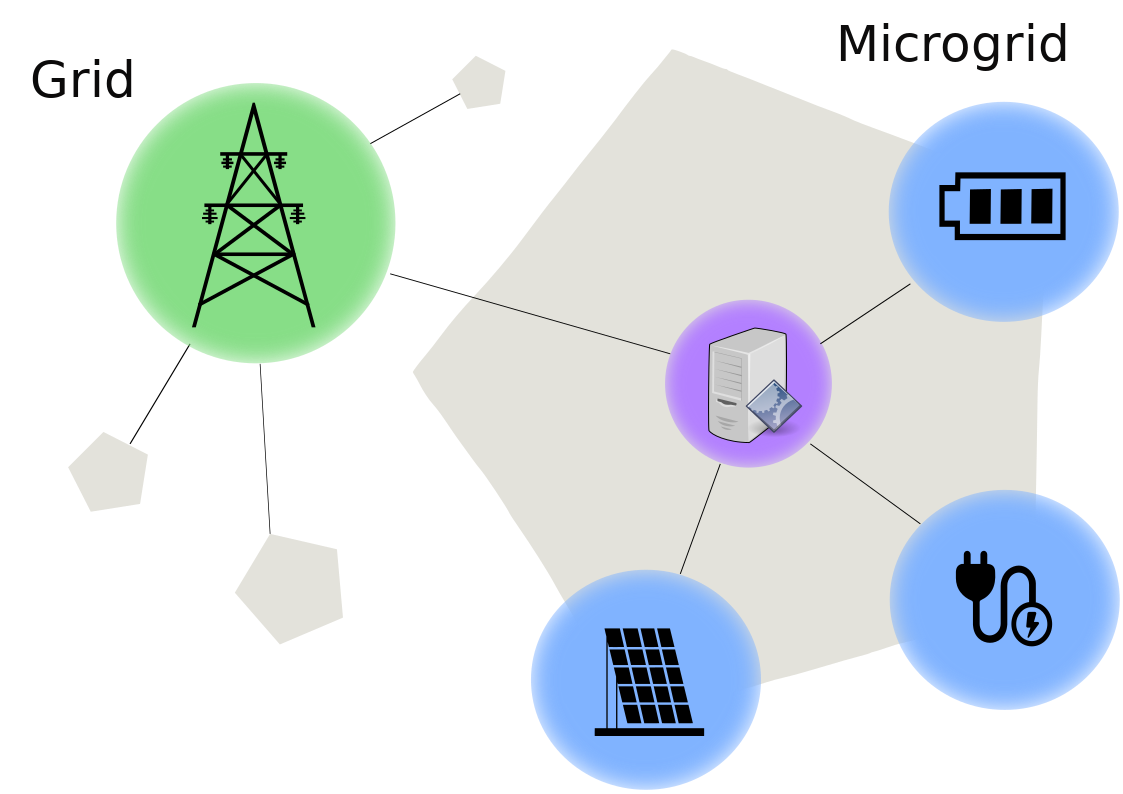}
	\caption{Schematic microgrid model}
	\label{fig:microgrid}
\end{figure}

We consider an electric microgrid composed of a photovoltaic power unit, an
electric load and an energy storage system. We assume that all components of
the microgrid share a single point of connection with the global grid. At this point, 
electric power can
be imported or exported so as to satisfy the electric power demand (total load)
at all times. We provide a schematic model of such a system
in Figure~\ref{fig:microgrid}. 

We now introduce some notation. 
We represent the mathematical control of a microgrid over a
finite number of discrete time steps $t \in \{0, 1, 2, \ldots, \horizon{-}1,
\horizon \}$, where
unit steps are spaced by~$\Delta_{t} = 15$ minutes. We denote a time interval 
between two decisions by~\( [t,t{+}1[ \), and not 
by~\( [t,t{+}1] \), to indicate that a decision is taken 
at the beginning of the time interval~\( [t,t{+}1[ \), 
and that a new one will be taken 
at the beginning of the time interval~\( [t{+}1,t{+}2[ \), 
and that these two consecutive intervals do not overlap.

\subsubsubsection{Storage dynamics}

\begin{subequations}
	The storage system is assumed
	to be a lithium-ion battery (or a container of aggregated batteries),
	characterized by the coefficients $\batt$ referring respectively to the
	battery's capacity (kWh), maximum load (kW), charge and discharge efficiency
	coefficients. The \emph{state of charge}, at the beginning of the time
	interval~$[t,t{+}1[$, is denoted by 
	\begin{equation}
	\soc_t\in \nc{0, 1} 
	\eqfinp \label{eq:state_of_charge}
	\end{equation}
	The \emph{decision~$\control_t$},
	taken at the beginning of every time interval~$[t,t{+}1[$,
	accounts for \emph{the energy charged} ($\control_t \geq 0$)
	or \emph{discharged} ($\control_t \leq 0$)
	during the time interval~$[t,t{+}1[$. 
	The dynamics of the state of charge is given by 
	\begin{align}
	\soc_{t{+}1}
	&=
	\dynamics\np{\soc_t, \control_t} 
	\eqsepv
	\forall t \in \na{0, \ldots,\horizon{-}1}
	\eqfinv
	\intertext{where the \emph{dynamics~$\dynamics$} is given by}
	\dynamics\np{\soc, \control}
	&=
	\soc + \frac{\rho_c}{c} \control^+ - \frac{1}{\rho_d c} \control^-
	\eqsepv \forall \np{\soc, \control} \in \nc{0, 1}\times \RR
	\eqfinv
	\end{align}
	\label{eq:dynamics}
\end{subequations}
with \( \control^+ = \max(0, u) \) and \( \control^- = \max(0, -u) \) .

\subsubsubsection{Constraints}

Constraints of the form
\begin{subequations}
	\begin{equation}
	\control_t \in \ConstraintSet\np{\soc_t}
	\eqsepv
	\forall t \in \na{0, \ldots,\horizon{-}1}
	\end{equation}
	restrict decisions $\control_t$ to the admissibility set
	(related to some of the battery parameters $\batt$ by means of the
	dynamics~$\dynamics$ in~\eqref{eq:dynamics}) 
	\begin{equation}
	\ConstraintSet\np{\soc} = \defset{\control \in \RR }%
	{\underline{\control} \leq \control \leq \overline{\control}
		\mtext{ and } 0 \leq \dynamics\np{\soc, \control} \leq 1}
	\eqfinv
	\label{eq:constraint_set}
	\end{equation}
	\label{eq:constraints}
\end{subequations}
where \( \overline{\control} = \overline{l} \times \Delta_{t} \) and
\( \underline{\control} = -\overline{l} \times \Delta_{t} \)
are the bounds on the energy that can be exchanged with the battery during a \noir{15-minutes} interval.

\subsubsubsection{Uncertainties}

For the purpose of defining the management costs,
we introduce the \emph{uncertainties}
\begin{equation}
\uncertain_t = \np{\pv_t, \demand_t} \in \RR^2
\eqsepv \forall  t \in \wspan
\label{eq:uncertainty}
\end{equation}
which represent a \noir{couple} of photovoltaic generation~$\pv_t$
and of energy demand~$d_t$.
As defined, uncertainties are \noir{exogenous} model variables.
\noir{When we turn to numerical experiments in~\S\ref{sec:experiments},
	sequences $\np{\uncertain_1, \ldots, \uncertain_{T}} \in \RR^{2\times T}$
	of uncertainties are obtained as
	samples from the historical observations of~\S\ref{Historical_observations}.}


\subsubsubsection{Costs}

\begin{subequations}
	We now turn to the management costs.
	The \emph{stage cost} during the time interval~$[t, t{+}1[$ is
	\begin{equation}
	\StageCost_{t}\np{\control_t, \uncertain_{t{+}1}} =
	p^+_{t} \cdot e_{t{+}1}^+ - p^-_{t} \cdot e_{t{+}1}^-
	\eqsepv \forall t \in \uspan
	\eqfinv
	\label{eq:cost}
	\end{equation}
	where
	\begin{equation}
	\ExchangeEnergy_{t{+}1} = \demand_{t{+}1} - \pv_{t{+}1} + \control_t
	\eqsepv \forall t \in \uspan
	\eqfinv
	\end{equation}
	is the \emph{energy exchanged with the grid} ---
	which, like the uncertainty~$\np{\pv_{t{+}1}, \demand_{t{+}1}}$,
	materializes at the end of the time interval~$[t,t{+}1[$,
	hence the index~$t{+}1$ ---
	and where $(p^+_{t}, p^-_{t})$ is the energy tariff
	(buying at price~$p^+_{t}$ and selling at price~$p^-_{t}$)
	applied during the time interval~$[t, t{+}1[$.
	
	Given a \noir{sequence $(\uncertain_1, \ldots, \uncertain_{T})$ of uncertainties} 
	and a sequence $(\control_0, \ldots, \control_{\horizon{-}1})$ of controls,
	we obtain the \emph{total operating cost}
	\begin{equation}
	\TotalCost\np{\control_0, \ldots, \control_{\horizon{-}1},
		\uncertain_1, \ldots, \uncertain_T}
	= \sum_{t=0}^{\horizon{-}1} \StageCost_t(\control_t, \uncertain_{t{+}1})
	\eqfinp
	\label{eq:total_cost}
	\end{equation}
	\label{eq:cost_structure}
\end{subequations}

In conclusion, we have introduced a dynamical system
with dynamics~\eqref{eq:dynamics},
constraints~\eqref{eq:constraints}
and a cost structure~\eqref{eq:cost_structure}.

\subsection{Microgrid controller}
\label{Microgrid_controller}

A microgrid controller is a mathematical device that, given some information at time~$t$,
yields a decision~$\control_t$.
We now detail the structure of the controllers that we will consider.
For this purpose, we introduce chronicles and management cost.

\subsubsubsection{Partial chronicles}

At the beginning of the time interval~$[t, t{+}1[$,
we \noir{may use} all the past observations and past forecasts
to make a decision~$\control_t$. For practical computational reasons,
we have chosen to restrict this information to
the \emph{partial observations}
\begin{subequations}
	\begin{equation}
	\np{ \uncertain_{t}, \uncertain_{t-1}, \ldots , \uncertain_{t-95} } =
	\Bp{\begin{array}{ll} & g_{t}, \ldots, g_{t-95} \\  & d_{t}, \ldots, d_{t-95} \end{array}}
	\in \RR^{2\times 96} \eqsepv \forall t \in \uspan\eqsepv
	\end{equation}
	of uncertainties~\eqref{eq:uncertainty} over the last 24 hours,
	and to the \emph{partial forecasts}
	\begin{equation}
	(\hat{\uncertain}_{t, t{+}1}, \ldots, \hat{\uncertain}_{t, t+96}) =
	\Bp{\begin{array}{ll} &\hat{g}_{t, t{+}1}, \ldots, \hat{g}_{t, t+96} \\  &\hat{d}_{t, t{+}1}, \ldots, \hat{d}_{t, t+96} \end{array}}
	\in \RR^{2\times 96} \eqsepv \forall t \in \uspan
	\eqfinv
	\end{equation}
	which represent a prediction 
	of the uncertainties~\eqref{eq:uncertainty} for the next 24 hours.
	Combined together, we obtain the \emph{partial observations-forecasts}
	\begin{align}
	\hist_t
	&=
	\Bp{\begin{array}{ll}
		&\uncertain_{t}, \uncertain_{t-1}, \ldots , \uncertain_{t-95} \\
		&\hat{\uncertain}_{t, t{+}1}, \ldots, \hat{\uncertain}_{t, t+96}
		\end{array}} \in \HIST \eqsepv \forall t \in \uspan \eqsepv
	\\
	\text{where } \quad              \HIST
	&=
	\RR^{2\times 96} \times \RR^{2\times 96}
	\eqfinv
	\end{align}
	\label{eq:partial_observations-forecasts}
	and, stacking them all over the whole time span, we obtain the \emph{partial chronicle}
	\begin{equation}
	\hist = \np{ \hist_0, \ldots, \hist_{\horizon{-}1}} \in \HIST^T
	\eqfinp
	\label{eq:partial_chronicle}
	\end{equation}
\end{subequations}

\subsubsubsection{Controller}

A \emph{controller~$\phi$} is a sequence of mappings
\begin{subequations}
	\begin{align}
	&\phi = \np{\phi_0, \ldots, \phi_{\horizon{-}1}}\eqsepv
	\intertext{where, for all $t \in \uspan$, we have 
		$\phi_t: \nc{0, 1} \times \HIST \to \RR$ with the constraints} 
	&\phi_t\np{\soc_t, \hist_t}
	\in \ConstraintSet\np{\soc_t}
	\eqsepv \forall \np{\soc_t, \hist_t} \in \nc{0, 1} \times \HIST
	\eqfinv
	\end{align}
	\label{eq:controller}
\end{subequations}
where the constraint set~\(  \ConstraintSet\np{\soc_t} \subset \RR \)
is defined in~\eqref{eq:constraint_set}.

\subsubsubsection{Management cost of a controller along a partial chronicle
	on a given site}

All dynamics~\eqref{eq:dynamics},
constraints~\eqref{eq:constraints}
and cost structure~\eqref{eq:cost_structure}
depend on parameters relative to a site.
Therefore, on a given site $\site \in \SITE$, 
we denote by~$\dynamics^{\site}$ the dynamics of the battery
in~\eqref{eq:dynamics}
and by \( \ConstraintSet^{\site} \) 
the set-valued mapping defining the
constraints as described in~\eqref{eq:constraint_set},
as they depend on the local parameters
\( (c^{\site}, \overline{l}^{\site}, \rho_c^{\site},\rho_d^{\site}) \).
We also denote by~$\StageCost^{\site}_t$ the stage cost
in~\eqref{eq:cost} --- as it depends on the energy tariff 
$(p^{+, \site}_{t}, p^{-,  \site}_{t})$ which could possibly be local 
--- and by~$\TotalCost^{\site}$ the total operating cost
in~\eqref{eq:total_cost},
as it depends on dynamics~$\dynamics^{\site}$,
set-valued mappings~\( \ConstraintSet^{\site} \) and
stage costs~$\StageCost^{\site}_t$.

Besides battery parameters and energy tariffs, sites differ from each other in
their historical data. For instance, the RMSE of the historical forecasts
(Figure~\ref{fig:RMSE}) varies across the pool of sites. Therefore,
controllers in~\eqref{eq:controller} may differ accordingly. 
This is why we denote by~$\phi^{\site}$ a controller for the 
site $\site \in \SITE$.

The application of a controller~$\phi^{\site}$ in~\eqref{eq:controller} 
along a partial chronicle~$\hist\in \HIST^T$
in~\eqref{eq:partial_chronicle} yields the \emph{management cost}
\begin{subequations}
	\begin{align}
	\Criterion^{\site}\np{\phi^{\site}, \hist}
	&=  \sum_{t=0}^{\horizon{-}1} \StageCost_t^{\site}(\control_t^{\site}, \uncertain_{t{+}1}) \eqsepv
	\label{eq:management_cost_controller_criterion}
	\intertext{where, for all $t \in \uspan$, the uncertainty
		$\uncertain_{t+1}$ is a component of $\hist_{t+1}$, and	
		the sequence $(\control_0^{\site}, \ldots, \control_{\horizon{-}1}^{\site})$ of controls
		is given by
	} 
	\soc_{0}^{\site}
	&=
	0 \eqfinv
	\label{eq:management_cost_controller_initial_state}
	\\
	\soc_{t{+}1}^{\site}
	&=
	\dynamics^{\site}\np{\soc_t^{\site}, \control_t^{\site}} \eqfinv
	\label{eq:management_cost_controller_dynamics}
	\\
	\control_t^{\site}
	&=
	\phi_t^{\site} \np{\soc_t^{\site}, \hist_t}\eqfinp
	\label{eq:management_cost_controller_feedback}
	\end{align}
	\label{eq:management_cost_controller}
\end{subequations}
The management cost~\( \Criterion^{\site}\np{\phi^{\site}, \hist} \)
will serve the assessment of the controller~$\phi^{\site}$
in~\S\ref{Designing_and_assessing_a_controller_on_a_given_site}.

\subsection{Designing and assessing a controller on a given site}
\label{Designing_and_assessing_a_controller_on_a_given_site}

We consider a given site $\site \in \SITE$.
We outline how to design and assess a controller~$\phi^{\site}$ as in~\eqref{eq:controller}.

\subsubsubsection{Data partitioning}

We have split the \noir{database} in periods of one week ranging from Monday 00:00 to
Sunday 23:45, each week containing thus $T=672=7 \times 24 \times 4 $ time
steps. With this, the \noir{database} is now organized as a subset~\(
\BaseData^{\site} \subset \HIST^{T} \) of chronicles, where $\HIST$ has been defined 
in~\eqref{eq:partial_observations-forecasts}.

Then, we partition the chronicles in the data set~\( \BaseData^{\site} \) in two disjoint subsets,
\( \CalibrationData^{\site} \) for calibration (training, in-sample)
and
\( \SimulationData^{\site} \) for simulation (testing, out-of-sample):
\begin{equation}
\CalibrationData^{\site} \cup \SimulationData^{\site} = \BaseData^{\site} 
\subset \HIST^{T}
\eqsepv
\CalibrationData^{\site} \cap \SimulationData^{\site} = \emptyset
\eqfinp
\label{eq:calibration_simulation_data_sets}
\end{equation}
For the EMSx benchmark, we select randomly 40\% of the weeks for simulation
and let the other 60\% be available for calibration.

\subsubsubsection{Calibration data}

The calibration data in \( \CalibrationData^{\site} \)
is available for the design of microgrid controllers as
in~\eqref{eq:controller}.
The design can result from any sort of technique (see examples in
\S\ref{Controllers_obtained_by_online_computation_methods}
and in
\S\ref{Controllers_obtained_by_offline_online_computation_methods}
below). 

\subsubsubsection{Simulation data}

On top of the weekly periods, every simulation chronicle in~\( \SimulationData^{\site} \)
is augmented with the data of the
Sunday before the period starts, following our
definition~\eqref{eq:partial_chronicle}.
Therefore, when simulating a microgrid
controller, 24~hours of past history data is always available to the decision-maker.
Simulation chronicles serve for testing only;
as such, they cannot be employed for the design of a controller.

\subsubsubsection{Parameters}
Additionally, we \noir{have} the battery parameters 
$\np{c^{\site}, \overline{l}^{\site}, \rho_c^{\site}, \rho_d^{\site}}$
(that have been adapted by Schneider Electric for the EMSx benchmark).
For computing the management cost in~\eqref{eq:management_cost_controller}, 
we use the energy tariff and time of use $(p^+_{t}, p^-_{t})$, in~\euro/kWh, from the French
electricity provider \'Electricit\'e de France (EDF);
it is the same for all sites.
At the beginning of any simulation, we assume the battery to be empty
(i.e. $\soc_0 = 0$), and we do not impose any final cost
(which is consistent with the expression~$\TotalCost^{\site}$ of the total
operating cost in~\eqref{eq:total_cost}).

\subsubsubsection{Lower bound for the management cost}

For any controller~$\phi^{\site}$ as in~\eqref{eq:controller}
and any partial chronicle~$\hist\in \HIST^T$ in~\eqref{eq:partial_chronicle},
the management cost~\( \Criterion^{\site}\np{\phi^{\site}, \hist} \)
in~\eqref{eq:management_cost_controller}
always has the so-called ``anticipative'' lower bound
\( \underline{\Criterion}^{\site}\np{\hist} \),
computed as the minimum of~\eqref{eq:management_cost_controller_criterion}
under the same constraints,
initial state~\eqref{eq:management_cost_controller_initial_state}
and dynamics~\eqref{eq:management_cost_controller_dynamics},
but where the last constraint
\( \control_t = \phi_t^{\site} \np{\soc_t, \hist_t} \)
in~\eqref{eq:management_cost_controller_feedback}
is now enlarged as
\( \control_t = \psi_t^{\site} \np{\soc_t, \hist} \),
for all $t \in \uspan$, for any 
\( \psi_t^{\site} : [0,1] \times \HIST^T \to \RR \). 
This gives a lower bound, because the minimization is done over 
``anticipative'' control laws
\( \psi_t^{\site} : [0,1] \times \HIST^T\to \RR \),
which encompass any controller~$\phi^{\site}$ as in~\eqref{eq:controller}.
Therefore, we easily get for any controller~$\phi^{\site}$, as in~\eqref{eq:controller}, that we have
\begin{equation}
\underline{\Criterion}^{\site}\np{\hist} \leq \Criterion^{\site}\np{\phi^{\site}, \hist} \eqsepv \forall \hist \in \SimulationData^{\site}\eqfinp
\label{eq:anticipative_upper_bound_score_site}
\end{equation}

\subsubsubsection{Performance score}

With the battery parameters and energy tariff, 
and a given site controller~$\phi^{\site}$ as in~\eqref{eq:controller},
the simulator in~\eqref{eq:management_cost_controller} yields as many
management costs~\( \Criterion^{\site}\np{\phi^{\site}, \hist} \)
as there are chronicles~$\hist$
in the simulation chronicles~$\SimulationData^{\site}$. 
Because the volume of energy production and consumption is variable from one site to
another, raw management costs~\( \Criterion^{\site}\np{\phi^{\site}, \hist} \) are not
suitable for a global performance analysis if we want to assess 
not only  a given controller, but a controller-design technique
(see examples in
\S\ref{Controllers_obtained_by_online_computation_methods}
and in
\S\ref{Controllers_obtained_by_offline_online_computation_methods},
and see \S\ref{Assessing_a_controller_design_technique} below).

Therefore, we shift the management cost of controller $\phi^{\site}$ 
over a chronicle $\hist \in \SimulationData^{\site}$ by subtracting the management 
cost of a dummy controller~$\phi^{\mathrm{d}}$
such that $\phi^{\mathrm{d}}_t = 0$, for $t \in \{0, \ldots, \horizon{-}1\}$. 
This dummy zero policy~$\phi^{\mathrm{d}}$ gives us a baseline operating cost,
the one that we would get when the microgrid is not 
equipped with an EMS and a battery. 
After the shift, we define 
the \emph{gain} of controller $\phi^{\site}$ by the following expression
\begin{equation}
\Gain^\site(\phi^\site) = \frac{1}{ \cardinal{ \SimulationData^{\site}} }
\sum_{ \hist \in \SimulationData^{\site} }
\bp{\Criterion^{\site}\np{\phi^{\mathrm{d}},\hist} -
	\Criterion^{\site}\np{\phi^{\site},\hist}} 
\eqfinv
\label{eq:gain}
\end{equation}
which expresses the average gain of introducing $\phi^\site$ in site $\site$,
over the baseline case of a dummy controller that does not use the battery. The
lower bound for the management cost
in~\eqref{eq:anticipative_upper_bound_score_site} provides a natural
\emph{upper bound for the gain}
\begin{subequations}
	\begin{align}
	&\overline{\Gain}^\site = \frac{1}{ \cardinal{ \SimulationData^{\site}} }
	\sum_{ \hist \in \SimulationData^{\site} }
	\bp{\Criterion^{\site}\np{\phi^{\mathrm{d}},\hist} -
		\underline{\Criterion}^{\site}\np{\hist}} \eqfinv
	\label{eq:upper_bound_for_the_gain}
	\\
	\text{where } \quad              
	&\Gain^\site(\phi^\site) \leq \overline{\Gain}^\site \eqfinp 
	\end{align}
\end{subequations}
In order to ease the performance analysis of a controller-design technique over an aggregated 
group of sites from $\SITE$, we scale the gain of $\phi^\site$ with the upper bound 
$\overline{\Gain}^\site$, and we define the \emph{performance score}
\begin{equation}
\Score^\site(\phi^\site) =
\frac{\Gain^\site(\phi^\site)}{\overline{\Gain}^\site} 
\eqfinp
\label{eq:score_site_controller}%
\end{equation}
The higher the score, the higher the gain allowed by the
controller~$\phi^{\site}$. A controller $\phi^\site$
that improves on the dummy controller~$\phi^{\mathrm{d}}$ gives a score in $[0, 1]$,
else the score is negative.

\subsection{Assessing a controller-design technique}
\label{Assessing_a_controller_design_technique}

\subsubsubsection{\noir{Controller-design} technique}

In addition to assessing a given controller,
we also aim at assessing a \emph{design technique for controllers}.
We provide examples of design techniques 
in
\S\ref{Controllers_obtained_by_online_computation_methods}
and in
\S\ref{Controllers_obtained_by_offline_online_computation_methods}.
In particular, when the same design technique is used across 
all sites, we will consider
a collection \( \sequence{\phi^{\site}}{\site \in \SITE} \)
of controllers derived from the application of a single design technique,
but adapted to each site $\site \in \SITE$.

\subsubsubsection{Performance score of a collection of controllers}

For a given collection \( \sequence{\phi^{\site}}{\site \in \SITE} \) 
of controllers, one per site, 
we average the performance score~\eqref{eq:score_site_controller} over sites, 
yielding the \emph{performance score} 
(of a collection \( \sequence{\phi^{\site}}{\site \in \SITE} \) of controllers)
\begin{equation}
\Score\bp{ \sequence{\phi^{\site}}{\site \in \SITE} } =
\frac{1}{ \cardinal{ \SITE } }
\sum_{ \site \in \SITE }\Score^{\site}\np{\phi^{\site}}
\eqfinp
\label{eq:score_controller_collection}
\end{equation}
When the controllers~$\phi^{\site}$ in the collection 
\( \sequence{\phi^{\site}}{\site \in \SITE} \)
have been designed by the same technique, 
the score \( \Score\bp{ \sequence{\phi^{\site}}{\site \in \SITE} } \) 
in~\eqref{eq:score_controller_collection} is a proxy 
to measure the performance of a controller-design technique over a large range of situations,
both in time of the year and in type of microgrid.
It permits an
immediate interpretation for practitioners interested in deploying 
a technique to design controllers on real microgrids. The higher the score, the higher the gain
allowed by the technique. 

\section{The EMSx benchmark software}
\label{sec:emsx_package}

After Schneider Electric's dataset presented in Sect.~\ref{sec:data},
and the mathematical framework presented in Sect.~\ref{sec:mathematical_framework},
the third component of the EMSx benchmark is the EMSx \emph{benchmark software},
that we have developed and implemented as a Julia \cite{bezanson2012julia} package named
\emsx\,.

\subsection{Requirements for the EMSx benchmark software}

Performing a simulation loop over the 70 sites $\site \in \SITE$, for computing 
the management costs~\eqref{eq:management_cost_controller}
of a collection \( \sequence{\phi^{\site}}{\site \in \SITE} \) of controllers
on the simulation chronicles of \( \sequence{\SimulationData^{\site}}{\site \in \SITE} \),
and finally obtaining the score \eqref{eq:score_controller_collection}, is a 
\noir{time-consuming} computing task. Indeed, our total pool of simulation chronicles
\( \sequence{\SimulationData^{\site}}{\site \in \SITE} \) gathers data from 2474 testing 
weeks, each of which requires $672$ ($7$ days in a week $\times$ $96$ decisions
in a day) calls to a controller, hence a total of 1.6 million calls for assessing a
single control technique.

We have developed a software, called \emsx, to ease the numerical assessment
of controllers in the context of the EMSx benchmark. Our first goal is to provide
an efficient and fast computing tool for running every simulation loop. 
We have chosen the Julia language, as it is a perfect candidate for processing 
the large amount of data made available in the EMSx benchmark dataset.
In particular, Julia makes it simple to distribute the simulation loop on several 
CPU cores, which enables us to release the \emsx\ package with parallel computing 
options. Moreover, Julia meets our second expectation, namely that our simulation 
software should be flexible enough to easily implement a large range of controllers 
as defined in~\S\ref{Microgrid_controller}.
We illustrate such flexibility in Sect.~\ref{sec:experiments},
where we outline the numerical results obtained with \emsx\ for various controller
design techniques.

\subsection{What the \emsx\ software package does}

Given a site $\site \in \SITE$, a controller~$\phi^{\site}$ in~\eqref{eq:controller}
and a partial chronicle~$\hist\in \HIST^T$ in~\eqref{eq:partial_chronicle}
(in practice, $\hist\in \SimulationData^{\site}$, the simulation chronicles
in~\eqref{eq:calibration_simulation_data_sets}),
the \emsx\ software returns the sequence of
states of charge of the battery, the stagewise costs,
and, above all, the management cost~\( \Criterion^{\site}\np{\phi^{\site}, \hist} \)
in~\eqref{eq:management_cost_controller}
and the corresponding computing time.
To this end, a \emsx\ user must provide
the implementation of her controller~$\phi^{\site}$ in Julia code following an API that we
briefly describe now.

Figure~\ref{fig:code} illustrates how the
dummy controller $\phi^{\mathrm{d}}_t = 0, \ t \in \{0, \ldots, \horizon{-}1\}$ can be implemented
and tested in a few lines of code. 
First, we define a new type for our controller, named
\texttt{DummyController}, as a subtype of a built-in \emsx\ type, named \texttt{AbstractController} (line 3).
Then, we implement the body of the \texttt{compute\_control} function (in Julia jargon, \emph{method})
for the new specialized \texttt{DummyController} type (line 5).
In the given example, the body is reduced to a simple zero expression,
as we have to implement a \emph{method} which always returns zero.
Eventually, we create an instance of a \texttt{DummyController} (line 8) and
launch the simulation over all simulation chronicles (line 10) passing the created instance as argument.
Battery parameters $\batt$ (referred to as metadata, line 13) and energy tariff
$(p^+, p^-)$ (line 12) can be changed by the user to run a custom
simulation. 

\begin{figure}[H]
	\centering
	\lstinputlisting{demo.jl}
	\caption{Example of the implementation and simulation of the dummy controller with the \emsx\ package}
	\label{fig:code}
\end{figure}

\noir{We can implement more complex examples by using
  the \texttt{information} object given in the input arguments of the
   \texttt{compute\_control} function.}
\noir{Figure~\ref{fig:information} displays the fields
of the \texttt{Information} type, an \emsx\ built-in type.
An instance of this type gives access to the running time step
$t \in \{0, \ldots, \horizon{-}1\}$ (line 2), to the state of charge
$\state_t$ in~\eqref{eq:state_of_charge} (line 3), and to the content
of the partial observations-forecasts $\hist_t$ 
in~\eqref{eq:partial_observations-forecasts} (line 4-7).
This online information allows us to define a controller $\phi$
as in~\eqref{eq:controller}. Besides, the \texttt{Information} type
has fields providing access to the energy price
$(p^+, p^-)$ (line 8), to the battery parameters $\batt$ (line 9)
and to the site reference $\site \in \SITE$ (line 10), which let us define
a specific controller $\phi^i$ for the site $\site \in \SITE$.
All-in-all, beyond our toy example of Figure~\ref{fig:code}, line 5-6,
the \texttt{information} object in the input argument 
of the \texttt{compute\_control} function 
enables the implementation of more sophisticated controllers,
such as the ones that we introduce and evaluate in~\S\ref{sec:experiments}.
}

\begin{figure}[H]
	\centering
	\lstinputlisting{information.jl}
	\caption{Detail of the \texttt{Information} type from the \emsx\ package}
	\label{fig:information}
\end{figure}

\noir{As an example of using the \texttt{Information} type, a controller that empties half of the battery would be implemented as 
in Figure~\ref{fig:code}, with the lines 5-6 replaced by
\lstinputlisting{demo2.jl}
}

The code of \emsx\ and more illustrative controller examples are
publicly available at \url{https://github.com/adrien-le-franc/EMSx.jl}.

\section{Numerical experiments}
\label{sec:experiments}

To illustrate how we can use the EMSx controller benchmark of
Sect.~\ref{sec:mathematical_framework},
we present several controllers and provide their scores.
These controllers~$\phi$ all share the same structure:
for any step~\(  t \in \uspan \),
the quantity \( \phi_t\np{\soc_t, \hist_t} \) in~\eqref{eq:controller} is not given by
an analytical formula, but as the solution of a \emph{reference optimization problem}.
\noir{Controller-design} techniques differ according to the nature
of this latter problem, which is solved at every step~\(  t \in \uspan \), that is,
\emph{online} (``on the fly'') as a function of the current available
quantities~\( \np{\soc_t, \hist_t} \), namely state of charge of the battery
and couples of observations-forecasts up to time~$t$ (see~\eqref{eq:partial_chronicle})
A comprehensive overview of such techniques and
algorithms is given in \cite{bertsekas1995dynamic}.

\noir{Before we start, we introduce some terminology.
First, we define a \emph{scenario} as a sequence 
$\na{\uncertain_t}_{t \in \mathcal{\horizon}}$ of uncertainties
\( \uncertain_t = \np{\pv_t, \demand_t} \) as in~\eqref{eq:uncertainty}, and
where $\mathcal{\horizon} \subseteq \na{1, \ldots, \horizon}$.
Second,
we stress the difference between \emph{open-loop} and \emph{closed-loop}
solutions and optimization problems. 
In an intertemporal optimization problem, one may look either for solutions
that are only functions of time (open-loop), or for solutions that are functions
of time \emph{and} of other variables that are available up to this time (closed-loop). Therefore, an optimization problem
is said to be open-loop (resp. closed-loop) when its solutions are open-loop (resp. closed-loop).
We refer to \cite[\S1.1.3]{carpentier2015stochastic}
for a discussion about the difference between open-loop and closed-loop controls.}

In~\S\ref{Controllers_obtained_by_online_computation_methods},
we present a class of so-called lookahead methods where the reference optimization problem
is multistage open-loop.
In contrast,
in~\S\ref{Controllers_obtained_by_offline_online_computation_methods},
we present a class of so-called cost-to-go methods 
where the reference optimization problem, solved online, 
is one-stage but depends on cost-to-go functions \noir{which, themselves, are the output 
of closed-loop optimization problems} which are computed \emph{offline}.
In~\S\ref{Results}, we comment the results obtained when applying the
controllers above on the EMSx controller benchmark.

\subsection{Controllers obtained by lookahead computation methods}
\label{Controllers_obtained_by_online_computation_methods}

In lookahead methods, one solves,
for every step~\(  t \in \uspan \),
a reference multistage (``looking ahead'' from the current step~$t$)
optimization problem which is
open-loop, be it deterministic (MPC) or stochastic (OLFC).
Stochastic (scenario-based) lookahead methods are limited by the
exponential growth of the computing time with respect to the number of
scenarios.

\subsubsection{Model Predictive Control}

The \emph{Model Predictive Control} (MPC) method
is one of the most famous lookahead techniques
\cite[Vol.1, \S6.1]{bertsekas1995dynamic}.
The MPC method mainly exploits the forecast data.
It yields a controller
\( \phi^{\MPC}=\np{ \phi_0^{\MPC}, \ldots, \phi_{\horizon{-}1}^{\MPC} } \)
by solving a sequence of \emph{multistage deterministic} optimization problems over a fixed\footnote{%
	When the horizon extends further than the period, we truncate the lookahead
	window to $\min(H, \horizon-t{+}1)$. \label{ft:truncation}} horizon~$H$ (in the
numerical application, $H = 96$)
\begin{equation}
\left\{
\begin{array}{ll}
&\begin{aligned}
& \control_t\opt \in \argmin_{\control_t}
\min_{ \np{\control_{t{+}1}, \ldots, \control_{t{+}H{-}1}} }
& & \sum_{s=t}^{t{+}H{-}1} \StageCost_{s}(\control_{s}, \hat{\uncertain}_{t, s+1}) \eqsepv \\
&&& \soc_{s+1} = f(\soc_{s}, \control_{s}) \eqsepv \forall s \in \{t, \ldots, t{+}H{-}1\} \eqfinv \\
&&& \control_{s} \in \ConstraintSet(\soc_{s}) \eqsepv \forall s \in \{t, \ldots, t{+}H{-}1\} \eqsepv
\end{aligned} \\[0.5cm]
& \phi_t^{\MPC}\np{\soc_t, \hist_t} = \control_t\opt
\eqfinv
\end{array}
\right.
\label{eq:MPC}
\end{equation}
where only the first value \( \control_t\opt \) of an
optimal sequence
\( \np{\control_{t}\opt, \control_{t{+}1}\opt, \ldots, \control\opt_{t{+}H{-}1}} \) is
kept.

When used in simulation, 
only the simulation (testing) chronicles
in~\( \SimulationData^{\site} \)
in the partition~\eqref{eq:calibration_simulation_data_sets} are used
in the reference optimization problem~\eqref{eq:MPC}, 
and not the calibration (training) chronicles in~\( \CalibrationData^{\site} \);
moreover, only a subvector (of available forecasted values) of the
whole vector~\eqref{eq:partial_observations-forecasts} of partial observations-forecasts
is used, namely
\( \np{ \hat{\uncertain}_{t, t{+}1}, \ldots,\hat{\uncertain}_{t, t{+}H{-}1} } \).

Due to the mathematical expressions for the
dynamics~\eqref{eq:dynamics}, the
constraints~\eqref{eq:constraints}
and the cost function~\eqref{eq:cost_structure},
the multistage deterministic optimization problem~\eqref{eq:MPC} formulates here as a \noir{linear program}.

\subsubsection{Open Loop Feedback Control}

The \emph{Open Loop Feedback Control} (OLFC) method
belongs to the family of online lookahead methods
and its approach is similar to~MPC, but for a stochastic component.
It yields a controller
\( \phi^{\OLFC}=\np{ \phi_0^{\OLFC}, \ldots, \phi_{\horizon{-}1}^{\OLFC} } \)
by a sequence of \emph{multistage open-loop stochastic} optimization problems over a fixed\textsuperscript{\ref{ft:truncation}} 
horizon~$H$ (in the numerical application, $H = 96$)
\begin{equation}
\left\{
\begin{array}{ll}
&\begin{aligned}
& \control\opt_t \in \argmin_{\control_t}
\min_{ \np{\control_{t{+}1}, \ldots, \control_{t{+}H{-}1}} }
& & \sum_{\scenario \in \SCENARIO} \pi^\scenario_t
\bgp{\sum_{s=t}^{t{+}H{-}1} \StageCost_{s}(\control_{s}, \uncertain^{\scenario}_{t, s+1})} \eqsepv \\
&&& \soc_{s+1} = f(\soc_{s}, \control_{s}) \eqsepv \forall s \in \{t, \ldots, t{+}H{-}1\} \eqsepv \\
&&& \control_{s} \in \ConstraintSet(\soc_{s}) \eqsepv \forall s \in \{t, \ldots, t{+}H{-}1\} \eqsepv
\end{aligned} \\[0.5cm]
& \phi_t^{\OLFC}\np{\soc_t, \hist_t} =
\control\opt_t
\eqfinp
\end{array}
\right.
\label{eq:OLFC}
\end{equation}
The reference optimization problem~\eqref{eq:OLFC} is open-loop because
the sequence of control variables \( \np{\control_{t{+}1}, \ldots, \control_{t{+}H{-}1}} \)
in the (second) $\min$ is not indexed by the scenario references~\( \scenario \in
\SCENARIO \). Only the first value \( \control_t\opt \) of an optimal sequence
\( \np{\control_{t}\opt, \control_{t{+}1}\opt, \ldots, \control\opt_{t{+}H{-}1}} \) is kept.

The reference optimization problem~\eqref{eq:OLFC} is \emph{stochastic} because of the scenarios
$(\uncertain^{\scenario}_{t, t{+}1},\) \( \ldots, \uncertain^{\scenario}_{t, t+96})_{\scenario \in \SCENARIO}$,
together with their probabilities $(\pi^\scenario_t)_{\scenario \in \SCENARIO}$.
These scenarios and their probabilities are built from two sources:
on the one hand, from the subvector~\( (\hat{\uncertain}_{t, t{+}1}, \ldots,
\hat{\uncertain}_{t, t{+}H{-}1} ) \) of available forecasted values given in the
whole vector~\eqref{eq:partial_observations-forecasts} of partial
observations-forecasts in the simulation (testing) chronicles
in~\( \SimulationData^{\site} \)
in~\eqref{eq:calibration_simulation_data_sets};
on the other hand, from partial
observations-forecasts in the calibration (training) chronicles
in~\( \CalibrationData^{\site} \)
in~\eqref{eq:calibration_simulation_data_sets},
from which we calibrate a scenario generation model.
In the forthcoming numerical experiments
in~\S\ref{Results}, we generate scenarios by modeling
\noir{the deviations from the net demand 24-hour forecast}  
as a Markov chain. 
We detail our scenario generation method in the Appendix
(\S\ref{Generation_of_scenarios_for_the_Open_Loop_Feedback_Control_algorithm}).
In numerical implementations, the
number of samples used varies between 10, 50 or 100 scenarios.

\subsection{Controllers obtained by cost-to-go computation methods}
\label{Controllers_obtained_by_offline_online_computation_methods}

In cost-to-go methods, one solves online,
for every step~\(  t \in \uspan \),
a reference single stage stochastic optimization problem,
which itself depends on cost-to-go functions, computed offline.
These functions are called cost-to-go because, ideally,
they map any state of the system to the optimal,
over strategies (closed-loop), expected future cost
from a given time step to the final horizon.
The Stochastic Dynamic Programming (SDP) method is the most famous of
cost-to-go computation techniques \cite{bertsekas1995dynamic}.
Cost-to-go methods are limited by the
exponential growth of the computing time with respect to the dimension of the
state space.
\noir{In SDP methods, we use random variables, defined on an abstract probability space
$\epro$, and designed by bold capital letters like~$\va{\Uncertain}$.}

\subsubsection{Stochastic Dynamic Programming}
\label{subsec:SDP}

 \noir{Whereas MPC and OLFC are based on an open-loop
reference optimization problem, 
\emph{Stochastic Dynamic Programming} is based on a closed-loop
reference optimization problem, like in~\eqref{eq:multistage_stochastic_optimization_problem}.}
\begin{itemize}
	\item In the \emph{offline phase of the \SDP algorithm},
	one computes a sequence
	of so-called \emph{value functions} $(V_t)_{t=0,1,\ldots,\horizon{-}1,T}$
	by the Bellman (or dynamic programming) equation,
	backward for \( t \in \uspan \),
	\begin{subequations}
		\begin{align}
		V_T\np{\soc}
		&= 0
		\eqsepv \nonumber 
		\\
		V_t\np{\soc}
		&=
		\min_{\control \in \ConstraintSet\np{\soc}}
		\sum_{\scenario \in \SCENARIO^\offline} \pi^{\offline,\scenario}_{t{+}1}
		\Bp{ \StageCost_{t}\np{\control, \uncertain_{t{+}1}^{\offline,\scenario}}
			+ V_{t{+}1}\bp{\dynamics\np{\soc, \control}} }
		\eqfinp
		\label{eq:SDP_offline}
		\end{align}
		\item In the \emph{online phase of the \SDP algorithm},
		one computes
		\begin{equation}
		\left\{
		\begin{array}{ll}
		&\begin{aligned}
		& \control_t\opt \in \argmin_{\control \in \ConstraintSet\np{\soc_t}}
		\sum_{\scenario \in \SCENARIO^\online} \pi^{\online,\scenario}_{t{+}1}
		\Bp{ \StageCost_{t}\np{\control, \uncertain_{t{+}1}^{\online,\scenario}}
			+ V_{t{+}1}\bp{\dynamics\np{\soc_t, \control}} }
		\end{aligned} \\[0.5cm]
		& \phi_t^{\SDPt}\np{\soc_t, \hist_t} = \control_t\opt
		\eqfinp
		\end{array}
		\right.
		\label{eq:SDP_online}
		\end{equation}
	\end{subequations}
\end{itemize}

The reference optimization problem~\eqref{eq:SDP_offline}--\eqref{eq:SDP_online}
is \emph{stochastic} because of the scenarios
\( \np{\uncertain^{\offline,\scenario}_{t{+}1}}_{\scenario \in \SCENARIO^\offline} \),
together with their probabilities $(\pi^{\offline,\scenario}_{t{+}1})_{\scenario \in
	\SCENARIO^\offline}$,
and
of the scenarios
\( \np{\uncertain^{\online,\scenario}_{t{+}1}}_{\scenario \in \SCENARIO^\online} \),
together with their probabilities $(\pi^{\online,\scenario}_{t{+}1})_{\scenario \in
	\SCENARIO^\online}$.

The scenarios indexed by \( \scenario \in \SCENARIO^\offline \),
and their probabilities, are built \emph{exclusively} from partial
observations-forecasts in the calibration (training) chronicles
in~\( \CalibrationData^{\site} \),
in~\eqref{eq:calibration_simulation_data_sets}.
The scenarios indexed by \( \scenario \in \SCENARIO^\online \),
and their probabilities, could additionally integrate partial
observations-forecasts in the simulation (testing) chronicles
in~\( \SimulationData^{\site} \)
in~\eqref{eq:calibration_simulation_data_sets}.
The reader will find details of our scenario generation method in
Appendix~\S\ref{Generation_of_scenarios_for_the_SDP_and_SDP-AR_algorithms}.

The reference optimization problem~\eqref{eq:SDP_offline}--\eqref{eq:SDP_online}
is closed-loop because, under proper assumptions,
it provides the optimal solution to the following
\emph{multistage stochastic optimization problem}
\noir{(where the minimum is over strategies~$\psi$ which depend both on time and on past uncertainties
as in~\eqref{eq:multistage_stochastic_optimization_problem_psi})}
\begin{subequations}
	\begin{align}
	& \min_{\psi}
	& & \EE
	\Bc{\sum_{t=0}^{\horizon{-}1} \StageCost_t(\va{U}_t, \va{\Uncertain}_{t{+}1})} \eqsepv \\
	&&& \va{X}_{t{+}1} = f(\va{X}_t, \va{U}_t) \eqsepv \forall t \in \{0, \ldots, \horizon{-}1\} \eqsepv \\
	&&& \va{X}_0 = x_0 \eqsepv \\
	&&& \va{\Control}_\stamp = \psi_t(\va{\Uncertain}_0, \ldots, \va{\Uncertain}_t) \eqsepv \forall t \in \{0, \ldots, \horizon{-}1\} \eqsepv 
	\label{eq:multistage_stochastic_optimization_problem_psi} \\
	&&& \va{\Control}_\stamp \in \ConstraintSet(\va{X}_t) \eqsepv \forall t \in \{0,
	\ldots, \horizon{-}1\}
	\eqfinp
	\end{align}
	\label{eq:multistage_stochastic_optimization_problem}
\end{subequations}
Problem~\eqref{eq:multistage_stochastic_optimization_problem} is optimally solved
by the Bellman equation \eqref{eq:SDP_offline} in the case where the random variables 
$(\va{\Uncertain}_0, \ldots, \va{\Uncertain}_T)$ (noise process) are
stagewise independent
\cite{bertsekas1995dynamic,carpentier2015stochastic}.

\subsubsection{\SDPAR}
\label{subsec:SDPAR}

To account for possible stagewise dependence in the
uncertainties (noise process in~\eqref{eq:multistage_stochastic_optimization_problem}), one can extend
the state space with the observations of the $i$-th net demand lags
$z_{t-i} = d_{t-i}-\pv_{t-i}$, $i=1, \ldots, k$.
This gives a new state~\( \tilde{\state}_t \) --- 
which is a compression of the information given by~\( (\state_t, \hist_t) \) --- 
and new elements of a model, like in~\S\ref{Microgrid_control_simulation},
with 
a new dynamics~\( \tilde{\dynamics}_t \),
new constraints~\( \tilde{\ConstraintSet}\), 
and a new stage cost~\( \tilde{L}_t \) as follows,
for all \( t \in \uspan \):
\begin{subequations}
	\begin{align*}
	\tilde{\state}_t
	& = ( \soc_t, z_t, \ldots, z_{t-k+1})
	\in [0,1] \times \RR^k
	\eqfinv
	\\
	\tilde{x}_{t{+}1}
	&=
	\tilde{\dynamics}_t\np{\tilde{x}, \control_t, \epsilon_{t{+}1}}
	\eqfinv
	\\
	\tilde{\dynamics}_t\np{\tilde{\soc}_t, \control_t, \epsilon_{t{+}1}}
	&=
	\Bgp{\begin{array}{ll}
		& \dynamics(\soc_t, \control_t) \\
		& \sum_{j=0,\ldots,k-1} \alpha_t^{j} z_{t-j} + \beta_t + \epsilon_{t{+}1} \\
		& z_{t}, \ldots, z_{t-k+2}
		\end{array}}
	\eqfinv
	\\
	\tilde{\ConstraintSet}(\tilde{x}_t) &= \ConstraintSet(x_t) \eqfinv \\
	\tilde{L}_t\np{\tilde{x}_t, \control_t, \epsilon_{t{+}1}}
	&=
	\StageCost_t\np{\control_t, \sum_{j=0,\ldots,k-1} \alpha_t^{j} z_{t-j}
		+ \beta_t + \epsilon_{t{+}1}}
	\eqfinp
	\end{align*}
\end{subequations}
The coefficients
\( \alpha_t^{0}, \ldots, \alpha_t^{k-1} \) and the additive terms~\( \beta_t \)
are the elements of an auto-regressive model of order~$k$ (noted AR($k$))
\begin{equation*}
\va{Z}_{t{+}1} = \sum_{j=0,\ldots,k-1} \alpha_t^{j} \va{Z}_{t-j} + \beta_t
+ \boldsymbol{\epsilon}_{t{+}1}
\eqsepv \forall t \in \uspan\eqsepv
\end{equation*}
for the net demand process $\va{Z}$.
When the error process~$\boldsymbol{\epsilon}$ is assumed to be stagewise
independent, the following algorithm provides an optimal solution to
the multistage stochastic optimization problem~\eqref{eq:multistage_stochastic_optimization_problem}.
\begin{subequations}
	\begin{itemize}
		\item In the \emph{offline phase of the \SDPAR algorithm},
		one computes a sequence of new value functions
		$(\tilde{V}_t)_{t=0,\ldots, T}$, backward for \( t \in \uspan \), by
		\begin{equation}
		\begin{aligned}
		\tilde{V}_T(\tilde{\state})
		&=
		0
		\eqsepv
		\\
		\tilde{V}_t(\tilde{\state})
		&=
		\min_{\control \in \tilde{\ConstraintSet}\np{\tilde{\soc}}}
		\sum_{\scenario \in \SCENARIO^\offline} \pi^{\offline,\scenario}_{t{+}1}
		\Bp{ \tilde{L}_{t}\np{\tilde{\state}, u, \epsilon_{t{+}1}^{\offline,\scenario}} +
			\tilde{V}_{t{+}1}\bp{\tilde{f}_{t}\np{\tilde{\state}, u,
					\epsilon_{t{+}1}^{\offline,\scenario}}} }\eqfinp
		\end{aligned}
		\end{equation}
		\label{eq:SDP-AR_offline} 
		\item In the \emph{online phase of the \SDPAR algorithm},
		one computes
		\begin{equation}
		\left\{
		\begin{array}{ll}
		&\begin{aligned}
		& \control_t\opt \in \argmin_{\control \in \tilde{\ConstraintSet}\np{\tilde{\soc}_t}}
		\!\! \sum_{\scenario \in \SCENARIO^\online} \!\! \pi^{\online,\scenario}_{t{+}1}
		\Bp{ \tilde{L}_{t}\np{\tilde{\state}_t, u, \epsilon_{t{+}1}^{\online,\scenario}} +
			\tilde{V}_{t{+}1}\bp{\tilde{f}_{t}\np{\tilde{\state}_t, u,
					\epsilon_{t{+}1}^{\online,\scenario}}} }
		\end{aligned} \\[0.5cm]
		& \phi_t^{\SDPARt}\np{\soc_t, \hist_t} = \control_t\opt
		\eqfinp
		\end{array}
		\right.
		\label{eq:SDP-AR_online} 
		\end{equation}
	\end{itemize}
\end{subequations}

\subsection{\noir{Using EMSx to compare controller-design techniques}}
\label{Results}

We now comment the results obtained when applying the 
controller-design techniques introduced in
\S\ref{Controllers_obtained_by_online_computation_methods}
and in
\S\ref{Controllers_obtained_by_offline_online_computation_methods}
to the EMSx benchmark. 
Numerical experiments were run on an Intel Core Processor
of 2.5 GHz with 22 GB RAM. We used the LP solver CPLEX 12.9 for MPC, OLFC and to
compute the lower bounds~\( \underline{\Criterion}^{\site}(\hist) \), 
\( \hist \in \SimulationData^{\site} \), \( \site \in \SITE,\)
in~\eqref{eq:anticipative_upper_bound_score_site}.
We have summarized our results in Table~\ref{tab:score}.

\noir{Our goal is to illustrate how EMSx enables a fine comparison of 
controller-design techniques.}
First, we focus on lookahead 
methods in~\S\ref{lookahead_performance}, second, we turn to cost-to-go methods 
in~\S\ref{cost_to_go_performance}, and finally, we give a comparative
analysis of the winner control techniques from both family of methods 
in~\S\ref{A_comparative_zoom_on_MPC_and_SDP-AR}. 
\noir{We conclude that, among our pool of candidate techniques, \SDPAR
stands out as the best microgrid control method on the EMSx benchmark.}

In what follows, we organize the discussion around three points:
i) performance scores \eqref{eq:score_controller_collection} 
(second column of Table~\ref{tab:score});
ii) (relative) gains\footnote{We recall that gains were defined relatively to the cost performance of a dummy controller.}~\eqref{eq:gain} disaggregated per site (Figures~\ref{fig:gain_ordered_by_rmse} 
and \ref{fig:gain_ordered_by_gap});
iii) CPU time performances.
Regarding CPU time, we report the \emph{online time} (fourth column of Table~\ref{tab:score}) as the average
computing time required to yield a single control~$\control_t$, and the
\emph{offline time} (third column of Table~\ref{tab:score}) as the average computing 
time that is to be spent prior to the simulation step for solving the online control problems defined
in \eqref{eq:MPC}, \eqref{eq:OLFC}, \eqref{eq:SDP_online} and \eqref{eq:SDP-AR_online}.

\begin{table}[htbp]
	\centering
	\begin{tabular}{l|lll}
		& \begin{tabular}{@{}c@{}}\textbf{Performance} \\ \textbf{score}\end{tabular} & \begin{tabular}{@{}c@{}}\textbf{Offline time} \\ (seconds)\end{tabular} & \begin{tabular}{@{}c@{}}\textbf{Online time} \\ (seconds)\end{tabular} \\ \hline
		\textbf{MPC} & 0.487 &  - & 9.82 $10^{-4}$  \\ \hline
		\textbf{OLFC-10} & 0.506 &  - & 1.14 $10^{-2}$  \\ \hline
		\textbf{OLFC-50} & 0.513 &  - & 8.62 $10^{-2}$  \\ \hline
		\textbf{OLFC-100} & 0.510 &  - & 1.87 $10^{-1}$  \\ \hhline{=|===}
		\textbf{SDP} & 0.691 &  2.67 & 3.09 $10^{-4}$ \\ \hline
		\textbf{SDP-AR(1)} & 0.794 &  38.1 & 4.44 $10^{-4}$  \\ \hline
		\textbf{SDP-AR(2)} & 0.795 &  468 & 5.55 $10^{-4}$  \\ \hline
		\textbf{Upper bound} & \textbf{1.0} &  - & -  \\ \hline
	\end{tabular}
	\caption{Scores (second column, \noir{the higher the better}) 
		\( \Score\bp{ \sequence{\phi^{\site}}{\site \in \SITE} } \)
		in~\eqref{eq:score_controller_collection}
		and time performances 
		(third and fourth column, \noir{the lower the better}) 
		for collections
		\( \sequence{\phi^{\site}}{\site \in \SITE} \) 
		of controllers $\phi^{\site}$ designed with techniques (first column) from
		\S\ref{Controllers_obtained_by_online_computation_methods}
		and
		\S\ref{Controllers_obtained_by_offline_online_computation_methods}
		on the EMSx benchmark. \noir{The symbol - indicates an irrelevant item}
		\label{tab:score}}
\end{table}

\subsubsection{Lookahead methods}
\label{lookahead_performance}

First, we examine performance scores (second column of Table~\ref{tab:score}).
Whereas MPC uses a single scenario (the forecast), 
OLFC uses multiple scenarios; 
this helps OLFC improving the average score of
MPC from 0.487 to 0.513. As expected, increasing the number of scenarios from
10 to 50 improves the score of OLFC. However, the performance remains stagnant
when pushing up to 100 scenarios, with a slight decrease of OLFC-100 to 0.510. 
We expect that the improvement between MPC and OLFC should be sensitive 
to the scenario generation method employed. With our method, the progress of OLFC 
is modest. 

Second, we turn to appraise the gains disaggregated per site.
In Figures~\ref{fig:gain_ordered_by_rmse} and \ref{fig:gain_ordered_by_gap},
we display the gains per site~\(\Gain^{\site}\np{\phi^{\site}} \) in~\eqref{eq:gain}.
We omitted OLFC-10 and OLFC-100 to improve readability, given that these methods perform
slightly worse than OLFC-50. We observe that OLFC-50 returns slightly higher gains than MPC for
58 out of 70 sites. We conclude that the stochastic approach of OLFC makes it a slightly more
accurate microgrid controller-design technique than MPC.

Third, we discuss the computing time.
Lookahead methods advantageously do not include an offline stage (third column of Table~\ref{tab:score}). 
However, except for MPC, they require a rather long online computing time 
(fourth column of Table~\ref{tab:score}), 
with an order of magnitude  between $10^{-4}$ and $10^{-1}$ seconds per call to the controller, for they call a LP solver at each time step.
We observe that OLFC is at least about ten times slower than MPC. As expected, the more we add sample scenarios, 
the longer the OLFC computing time. We observe that improving the amount of scenarios from 50
to 100 doubles the online time. Even though the online time of OLFC-100 is 
still reasonable for a field implementation of a microgrid controller, it took us 
about 46 hours to run the simulation over the 2474 simulation weeks. Despite the 
much longer computing time, OLFC-100 did not return higher gains than OLFC-50. For 
this reason, we did not consider more than 100 generated scenarios.

\subsubsection{Cost-to-go methods}
\label{cost_to_go_performance}

\begin{figure}[htbp]
	\centering
	\includegraphics[width=\textwidth]{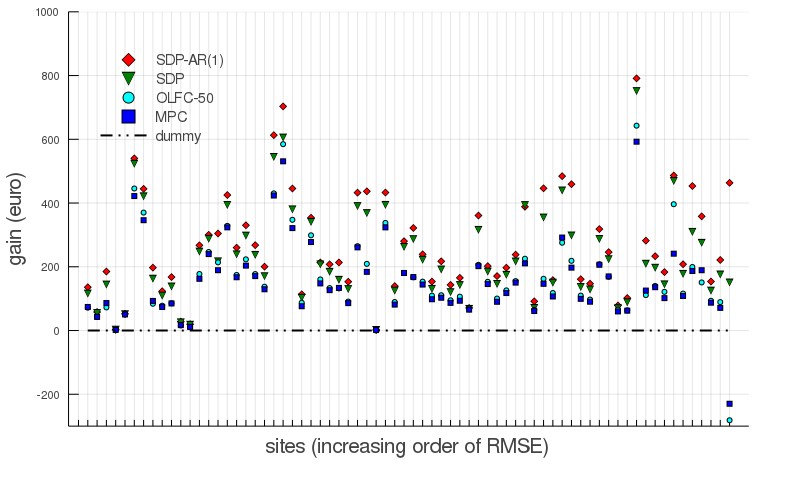}
	\caption{Gain $\Gain^{\site}(\phi^{\site})$  ($Y$-axis)
		\noir{in~\eqref{eq:gain}}
		per sites $\site \in \SITE$  ($Y$-axis) of MPC, OLFC-50,
		\SDP and \SDPAR(1); sites are
		ranked in increasing order of RMSE along the $X$-axis
		\label{fig:gain_ordered_by_rmse} }
\end{figure}

\begin{figure}[htbp]
	\centering
	\includegraphics[width=\textwidth]{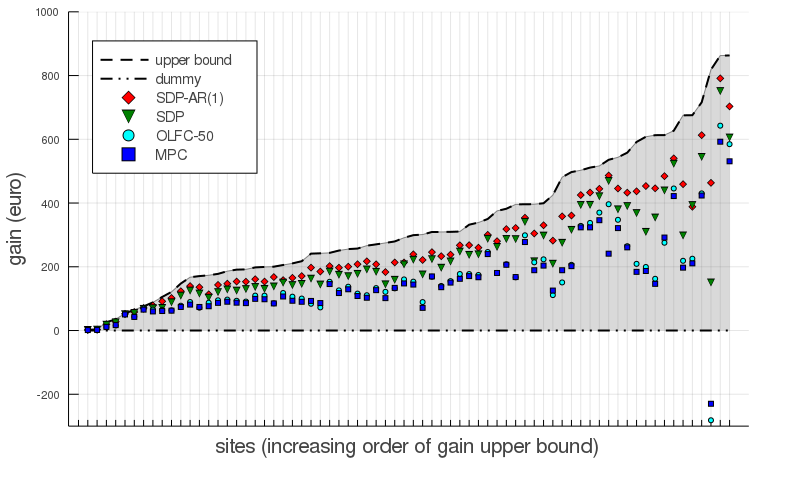}
	\caption{\noir{
		Gains ($Y$-axis) per sites $\site \in \SITE$ ($X$-axis):
		the upper bound (dashed line)
		is the quantity
		$\overline{\Gain}^\site$
		\noir{in~\eqref{eq:upper_bound_for_the_gain}};
		the horizontal dotted and dashed line represents the null gain of
		the dummy controller;
		the other symbols represent 
		the quantities
		$\Gain^{\site}(\phi^{\site})$
		in~\eqref{eq:gain}
		corresponding to the four methods
		\SDPAR(1), \SDP, OLFC-50, MPC}; sites are
	ranked in increasing order of gain $\overline{\Gain}^\site$ along the $X$-axis} 
	\label{fig:gain_ordered_by_gap}
\end{figure}

First, we examine performance scores (second column of Table~\ref{tab:score}).
We see that a plain \SDP controller yields scores jumping to 0.691.
The results of \SDP can be
improved up to 0.794 in the \SDPAR formulation. However, We observe that the gain from
extending the lag of an AR(1) model to an AR(2) model is almost null. Therefore,
we do not report the gains of \SDPAR(2) in Figures~\ref{fig:gain_ordered_by_rmse} 
and \ref{fig:gain_ordered_by_gap}. 

Second, we turn to appraise the gains disaggregated per site.
Looking closer at the per site performances, \SDPAR(1) outperforms
\SDP for 69 of the 70 sites. As suggested by its performance score, 
Figure~\ref{fig:gain_ordered_by_gap} reveals that the gains allowed 
by \SDPAR(1) are close to the upper bounds $\overline{\Gain}^\site, \site \in \SITE$
in~\eqref{eq:upper_bound_for_the_gain}.

Third, we discuss the computing time.
Cost-to-go methods display fast online times (fourth column of Table~\ref{tab:score}),
with an order of magnitude between  $10^{-3}$ and $10^{-4}$ seconds per call to the controller. 
However, cost-to-go methods require offline CPU time 
(third column of Table~\ref{tab:score}) for computing value functions.
The complexity of \SDP is well known for growing exponentially with the state space, 
which is well illustrated by our results: from \SDP to \SDPAR(1) to \SDPAR(2), we
add one dimension to the state space at each improvement of the method, which
multiplies by a factor of 10 the offline time. Given the low improvement of gain 
from \SDPAR(1) to \SDPAR(2), we find the offline time of the latter method dissuasive.
Finally, \SDPAR(1) appear as the best trade-off between computing time and cost performance.

\subsubsection{A comparison between lookahead and cost-to-go methods}
\label{A_comparative_zoom_on_MPC_and_SDP-AR}

We now discuss the comparative performances of OLFC-50 and \SDPAR(1). Both techniques
represent the best candidate of its family of method. 

\subsubsubsection{Cost performance}
First, we examine performance scores (second column of Table~\ref{tab:score}).
The main observation is that the
performance score of \SDPAR(1) is more than 50\% higher than the one
of OLFC-50 (second column of Table~\ref{tab:score}). 

Second, we turn to appraise the gains disaggregated per site.
The observed aggregated dominance 
is confirmed when looking closer at the per site performances in
Figures~\ref{fig:gain_ordered_by_rmse} and \ref{fig:gain_ordered_by_gap}.
Indeed, we see that \SDPAR(1) outperforms OLFC-50 for 68 of the total
pool of 70 sites. Even more, these two figures reveal that lookahead methods 
lag behind cost-to-go ones for almost all sites, and that the gap can be 
significant on a few outlying sites. The underperformance of OLFC-50 
on these sites explains the score gaps of Table~\ref{tab:score}.

Third, we discuss the relationship between the cost performances of both techniques
and the predictability of the sites, measured by the RMSE \eqref{eq:RMSE}. 
For this purpose, we comment and detail Figure~\ref{fig:score_vs_rmse}.
We immediately observe that, for \SDPAR(1), 
there is no strong link between cost performance and predictability.
Regarding OLFC-50, a statistical analysis reveals that the correlation 
between the RMSE value and the OLFC-50 score is moderate, 
with a Pearson correlation coefficient of -0.54
(the lower the RMSE, the higher the score).
For most of the sites, the RMSE explains very little of the performance
of OLFC-50, to the point that OLFC-50 performs poorly on some sites 
with low RMSE values (and thus high forecast accuracy). 
Now, we focus on a few number of sites that appear as outliers: 
Site~33 and Site~59 (the two circles at the far right),
and Site~48 (the circle at the top left).
OLFC-50 achieves its highest scores 0.965
on Site~33 and 0.910 on Site~59, which both have very regular and easily
predictable load profiles,
and its worse score -0.340 on Site~48, the least predictable 
site of the pool, with a RMSE of 0.17.
Notice that this negative score of -0.340 means that OLFC-50 performs worse
on average than the dummy controller~$\phi^{\mathrm{d}}$.
In contrast, we observe that \SDPAR(1) scores 0.566 on Site~48, which highlights
the robustness of the cost-to-go methods for the management of a microgrid 
in a highly unpredictable context.

\subsubsubsection{Computing time performance}
Second, we discuss the computing time.
On the one hand, regarding offline time, 
\SDPAR(1) cannot do better than OLFC-50, obviously, but 
the average offline time of \SDPAR(1) (38.1 seconds)
is reasonable for a field implementation of a microgrid controller. 
On the other hand, regarding online time, 
OLFC-50 is not as good as \SDPAR(1) --- 
with an average online time of OLFC-50 about 100 times longer than the
one of \SDPAR(1) --- but its order of magnitude ($10^{-2}$ seconds) remains acceptable. 
Thus, all in all, computing time is not a discriminating point 
between the two techniques.

\begin{figure}[htbp]
	\centering
	\includegraphics[width=0.9\textwidth]{score_olfc_vs_sdp_ar.png}
	\caption{Performance score $\Score^{\site}(\phi^{\site})$ per sites $\site \in \SITE$ 
		for the \SDPAR(1) and OLFC-50 computation methods ($X$-axis) \emph{versus} RMSE ($Y$-axis);
		the vertical line recalls the baseline score positioning of a dummy controller} 
	\label{fig:score_vs_rmse}
\end{figure}

\section{Conclusion}

We have introduced EMSx, an Energy Management
System benchmark to compare electric microgrid controllers,
hence controller-design techniques. 
EMSx is made of three key components. 
The dataset provided by Schneider Electric ensures a diverse pool of realistic
microgrids with photovoltaic power integration.
The mathematical framework is explicit.
The simulation code is accessible in the \emsx\ package, 
designed to welcome various sorts of control algorithms. 
All components of the benchmark are publicly
available, so that other researchers willing to test controllers on EMSx may
reproduce experiments easily.

Regarding our numerical results, we observe a gap
between cost-to-go methods and lookahead methods. The \SDPAR(1) controller-design technique
stands out as the best trade-off between cost optimality and computing time
performance. However, there is a range of possible improvements to explore.
Among interesting directions, other scenario generation
techniques could be tested to see how does the OLFC controller react. Beyond
plain score improvements, enriching contributions could arise from the reduction
of the computing time, 
or from changing the performance metrics (for instance with the use of risk measures),
and from further comparative analysis of methods, especially to better explain
the gap between lookahead and cost-to-go methods. We are also looking forward to
controllers inspired from other research fields than multistage deterministic or
stochastic optimization, including 
heuristics, robust optimization and reinforcement learning.

\begin{acknowledgements}
	We thank Efficacity and Schneider Electric for the PhD funding of Adrien Le
	Franc. Additionally, we are grateful for the feedbacks and data supply from
	Peter Pflaum and Claude Le Pape (Schneider Electric) and we thank our colleague
	Tristan Rigaut (Efficacity) for insightful tips about the Julia
        language. \noir{We thank the Guest Editor and the Reviewers for their
          insightful comments that helped improve the manuscript.}
        
\end{acknowledgements}

\appendix
\section{Scenario generation}
\label{Appendix}

\subsection{Generation of scenarios for the Open Loop Feedback Control
	algorithm}
\label{Generation_of_scenarios_for_the_Open_Loop_Feedback_Control_algorithm}

\begin{figure}[htbp]
	\centering
	\includegraphics[width=0.75\textwidth]{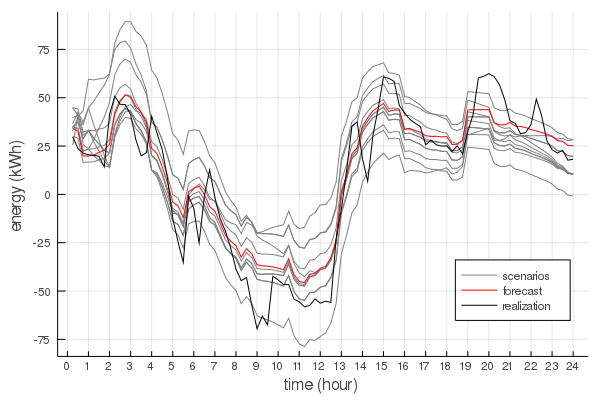}
	\caption{Example of scenarios generated from a 24 hours net demand forecast}
	\label{fig:scenarios}
\end{figure}

Our scenario generation method is inspired by \cite{staid2017generating,
	woodruff2018constructing}, which address such generation for so-called day-ahead
energy management problems. We use a simplified model for tractable generation adapted to
dynamical control. Since uncertainties $\uncertain_{t} = (\pv_t, \demand_t)$ are
directly plugged in the cost \eqref{eq:cost} as the net demand $z_{t} =
\demand_t - \pv_t$, we compress the generation of uncertainty scenarios
$\np{\uncertain^{\scenario}_{t, t{+}1}, \ldots, \uncertain^{\scenario}_{t, t+96}}
\in \RR^{2\times 96}$ in \eqref{eq:OLFC} to the generation of net demand scenarios
$\np{z^{\scenario}_{t, t{+}1}, \ldots, z^{\scenario}_{t, t+96}} \in \RR^{96}$.
Following the original methods, we choose day part separators to reduce the 96
dimensional vector to a few skeleton points for scenario  
construction (intermediate values are linearly interpolated). We choose to concentrate on the forecast error at $t$+15 minutes,
$t$+1 hour, $t$+2 hours, $t$+4 hours, $t$+12 hours and $t$+24 hours,
so that our model only samples values of $z^{\scenario}_{t, t{+}j}$
for $j \in \na{1, 4, 8, 16, 48, 96}$.
We select 
10 relevant values of the net demand error at each separator $j \in \na{1, 4, 8, 16, 48, 96}$ by applying the
$K$-means algorithm \cite[\S13]{friedman2001elements} on the historical error $z_{t+j}^{\site} - \hat{z}_{t, t+j}^{\site}$, combining the historical observations \eqref{eq:historical_observations} and forecasts \eqref{eq:historical_forecasts} of the calibration data of the Site $\site \in \SITE$ considered.
Then, we compute 10$\times$10 transition matrices for the error between
consecutive separators. With this model, we are able to sample net
demand scenarios $\np{z^{\scenario}_{t, t{+}1}, \ldots, z^{\scenario}_{t, t+96}}$ given a single value forecast and to compute their
probabilities~$\pi^{\scenario}_t$. Separate probability distributions of the initial error at $t$+1
are calibrated depending on the time of the day, and on whether the initial time step $t$
corresponds to a week day or a weekend day. We alleviate computing costs by
reusing transition matrices regardless of the initial value of~$t$. However we
calibrate separate matrices for week days and weekend
days. Figure~\ref{fig:scenarios} provides examples of scenarios generated with
our method.

\subsection{Generation of scenarios for the \SDP and \SDPAR algorithms}
\label{Generation_of_scenarios_for_the_SDP_and_SDP-AR_algorithms}

For \SDP methods, we choose to discretize each dimension of the state space in 10~values,
whereas the control space is restricted to 20~values and the noise space to
10~values. Since uncertainties $\uncertain_{t} = (\pv_t, \demand_t)$ are
directly plugged in the cost \eqref{eq:cost} as the net demand $z_{t} =
\demand_t - \pv_t$, we compress the calibration of the distributions of
$\np{\va{\Uncertain}_1, \ldots, \va{\Uncertain}_T}$ to the calibration of the
distributions of $\np{\va{Z}_1, \ldots, \va{Z}_T}$. We use the $K$-means
algorithm to fit discrete probabilities~$\pi^{\offline,\scenario}_{t{+}1}, \pi^{\online,\scenario}_{t{+}1}$ on the historical
observations $\np{\pv^{\site}, \demand^{\site}}$
\eqref{eq:historical_observations} of the calibration data of the Site $\site
\in \SITE$ considered. 
These discrete distributions serve the computation of expectations
in~\eqref{eq:SDP_offline}--\eqref{eq:SDP_online}. While we could leverage the
data available on the fly in the \emph{online phase}, we use the same
probability distributions in the \emph{offline phase} and in the \emph{online phase}. Separate
distributions (of~$\va{Z}_t$) are calibrated depending on the time of the
day and on whether the time step $t{+}1$ corresponds to a week day or a weekend day. We
compute one value function per site for the horizon of one week. In the \SDPAR
formulation, we calibrate the AR($k$) models using least square regression and
calibrate distributions of the residual error ($\boldsymbol{\epsilon}_t$) with the same approach as
for $\va{Z}_t$.

\bibliographystyle{plain}
\bibliographystyle{plain}

\end{document}